\documentclass[2p]{elsarticle}
\usepackage{graphicx}
\usepackage{amssymb}
\usepackage{amsmath}
\newtheorem{thm}{Theorem}
\newtheorem{lem}[thm]{Lemma}
\newtheorem{rmk}{Remark}
\newtheorem{ass}{Assumption}
\newproof{pf}{Proof}

\begin{document}

\begin{frontmatter}

\title{Finite Volume Difference Scheme for a Degenerate Parabolic Equation in the Zero-Coupon Bond Pricing\tnoteref{t2}}


\author{T. Chernogorova}
\ead{chernogorova@fmi.uni-sofia.bg}
\author{R. Valkov\corref{cor1}}
\ead{rvalkov@fmi.uni-sofia.bg} \cortext[cor1]{Corresponding
author}
\address{Faculty of Mathematics and Informatics, University of Sofia, 1000 Sofia, Bulgaria}

\begin{abstract} 

In this paper we solve numerically a {\emph{degenerate}}
parabolic equation with {\emph{dynamical}} boundary
conditions of zero-coupon bond pricing. First, we discuss
some properties of  the differential equation. Then, starting
from the divergent form of the equation we implement the finite-volume method
of S. Wang \cite{W} to discretize the differential problem.
We show that the system matrix of the discretization scheme
is a $M$-matrix, so that the discretization is {\emph{monotone}}.
This provides the non-negativity of the price with respect to time
if the initial distribution is nonnegative. Numerical experiments demonstrate the efficiency of our difference scheme
near the ends of the interval where the degeneration occurs.

\end{abstract}
\end{frontmatter}

\section{Introduction}

Since the Black-Scholes models rely on stochastic differential
equations, option pricing rapidly became an attractive topic for
specialists in the theory of probability and stochastic methods
were developed first for practical applications, along with
analytical closed formulas. But soon, with the rapidly growing
complexity of the financial products, other numerical solutions
became attractive [1,2,6,12,15-19].

There is a large and ever-going number of different interest rate derivative products now, for instance
bonds, bonds options, interest rate caps, swap options, etc. Bonds in general
carry coupons, but there also exists a special kind of bond without coupons which
is called zero coupon bond (ZCB). A ZCB is purchased today a certain price,
while at maturity the bond is redeemed
for a fixed price. By a similar way to the derivation of the Black-Sholes equation, the problem of ZCB pricing can be reduced to a partial differential equation (see [5,13]).

The present paper deals with a degenerate parabolic equation of
zero-coupon bond pricing [5,13]. Since our equation (see (1), (2), (3)) in the next section
becomes {\it {degenerate}} at the boundary of the domain,
classical finite difference methods may fail to give accurate
approximations near the boundary. An effective method that resolves the
singularity is proposed by S. Wang \cite{W} for the Black-Sholes
equation. The method is based on a finite volume formulation of
the problem coupled with a fitted local approximation to the
solution and an implicit time-stepping technique. The local
approximation is determined by a set of two-point boundary value
problems defined on the element edges. This fitting technique is
based on the idea proposed by Allen and Southwell [8,10] for
convection-diffusion equations and has been extended to one and
multidimensional problems by several authors [7,8,10].

This paper is organized as follows. Our model problem is presented
in Section 2, where we discuss our basic assumptions and some
properties of the solution. The discretization method is developed
in Section 3. Section 4 is devoted to the time discretization. We
show that the system matrix is a $M$-matrix, so that the
discretization is monotone. In this case the maximum principle is
satisfied and thus the discrete solution is non-negative. Numerical experiments
show higher accuracy of our scheme in comparison with
other known scheme near the degeneracy.
We observe and emphasize the fact that in the proposed method, we do not need to refine the
mesh near the boundary (degeneration).


\section{The continuous problem }

Suppose that the short term {\emph{interest rate, the spot rate,}} follows a random walk
\[
dr=\theta(r)dt+w(r)dz,
\]
where $z(t)$ is the Brownian motion. Since the {\emph {spot rate}}, in practice,
is never greater than a certain number, which is assumed $R$, and never
less than or equal to zero, we suppose that  $r \in [0,R]$.

\begin{ass} $\theta(r)$ is a Lipschitz function, which satisfies

\begin{equation}\label{1}
\theta(0) \geq 0, \; \;\theta(R)\leq 0.
\end{equation}

\end{ass}

\begin{ass} $w(r)$ is a non-negative and smooth bounded function, which satisfies

\begin{equation}\label{2}
w(0)=w(R)=0, \; \; w(r)>0, \; \; r\in (0,R).
\end{equation}

\end{ass}

By the well-known delta hedging technique, the ZCB premium
$P=P(r,t;T)$ satisfies the following backward parabolic equation (see [13]):
\begin{equation}\label{3}
\frac{\partial P}{\partial t}+\frac{w^{2}(r)}{2}\frac{\partial^{2} P}{\partial r^{2}}+
(\theta(r)+\lambda(t)w(r))\frac{\partial P}{\partial r}-rP = 0,\;\;
(r,t) \in Q \equiv [0,R] \times [0,T),
\end{equation}
\begin{equation}\label{4}
P(r,T)=Z,
\end{equation}
where $T$ is the {\emph{maturity}}, $Z$ is a fixed constant.
Function $\lambda(t)$ in (\ref{3}) is called the {\emph{market price risk.}}
For the given functions $\theta$, $w$ and $\lambda$, the problem of ZCB
pricing consists of the determination of the solution
$P(r,t)$ from equation (\ref{3}), which is often referred to as a direct problem.

Being different from the classical parabolic equations in which the
principal coefficient is assumed to be {\emph{strictly positive}},
the parabolic equation (\ref{3})
belongs to the second order differential equations with
non-negative characteristic form. The main character of such
kinds of equations is {\emph{degeneracy.}} It can be
easily seen that at
$r=0$ and $r=R$, equation (\ref{3}) degenerates into a hyperbolic equation with
positive and negative characteristics respectively
\begin{equation}\label{5}
\frac{\partial P}{\partial t}+\theta(0)\frac{\partial P}{\partial r}=0,
\end{equation}
\begin{equation}\label{6}
\frac{\partial P}{\partial t}+\theta(R)\frac{\partial P}{\partial r} =RP.
\end{equation}

By the Fichera's theory (see [9]) for degenerate parabolic equations, we
have that at the degenerate boundaries $r=0$ and $r=R$, the boundary conditions should not
be given. Therefore, the maturity  data $P(r,T)$ determines the solution $P(r,t)$ of
problem (\ref{3}), (\ref{4}) {\emph{uniquely}}.

First, we make the change of variable ${\widehat {t}}=T-t,$ and let
${\widehat{\lambda}}(t)=\lambda(T-\tau)$. Then, coming back to $t$, the function $P$ satisfies
the following parabolic equation
\begin{equation}\label{7}
\frac{\partial P}{\partial t}-\frac{w^{2}(r)}{2}\frac{\partial^{2}P}{\partial r^{2}}-
(\theta(r)+\lambda(t)w(r))\frac{\partial P}{\partial r}+rP = 0, \; \;(r,t) \in Q
\end{equation}
with  initial condition
\begin{equation}\label{8}
P(r,0)=P_{0}(r).
\end{equation}
Let us note that for the concrete model (\ref{3}),(\ref{4}) we consider $P_{0}(r)= Z$.

If the functions $\theta$, $w$ satisfy the Assumptions 1,2 and the
initial data $P_{0} (r) $ is a continuous function then there exists a classical solution
($P$ has continuous first derivative with respect to $t$
and second derivative with respect to $r$ up to the boundary
$\partial Q $ and satisfies equation (\ref{7}), see [5,9]) of the problem (\ref{7}), (\ref{8}).
Further, in Section 4, we show that our difference scheme satisfies a discrete analogue of the following maximum principle, see [5,9]:
\begin{lem}
Let Assumptions 1, 2 hold. Then
\[
0 \leq P(r,t) \leq P_{0}(r).
\]
\end{lem}

The Dirichlet problem on the domain $(0,X) \times (0,T), \; 0< X < \infty $
for the Black-Scholes equation [3], studied by Song Wang
[15, equation
(2.9a)], has the form (\ref{3}) with coefficients:
\[
\frac{1}{2}\sigma^{2}r^{2}\; \; \mbox{at} \; \;
\frac{\partial^{2} P}{\partial r^{2}}\; \mbox{and} \;
(d(t)-D(x,t))r \; \; \mbox{at} \; \;\frac{\partial P}{\partial r}.
\]
In [15] $P$ denotes the value of a European call or put option,
$\sigma= const>0$ denotes the volatility of the asset,
the interest rate are denoted by $r$ and $D$ are the dividends. It is assumed that $r>D$.
Following this line we will assume the further
specifications on $w(r)$ and $\theta(r)$:
\begin{equation}\label{9}
w(r)=r(R-r)w_{0}(r),
\end{equation}
where $w_{0}(r)\geq w_{0}= const>0$ is smooth function and
\begin{subequations}\label{10}
\begin{align}
\theta(r)&=r(R-r)\theta_{0}(r),\;\;\; \theta_{0}(0)\not=0,\; \theta_{0}(R)\not=0,\label{10a}\\
\theta(r)&=r\theta_{0}(r),\; \; \; \theta_{0}(R)<0,\label{10b}\\
\theta(r)&=(R-r)\theta_{0}(r),\; \; \; \theta_{0}(0)> 0,\label{10c}\\
\theta(r)&=\theta_{0}(r),\;\; \; \theta_{0}(0)>0,\;\theta_{0}(R)<0.\label{10d}
\end{align}
\end{subequations}

Let us note that similar to (\ref{7}) degenerate parabolic equations with
coefficients of type (\ref{9}), \eqref{10a}-\eqref{10d} are obtained by introducing new variables to
transform the problem posed on infinite to finite interval for discretely sampled Asian options
[1,19]. Also, see the models in [14].

Further, we will work with the following fully-conservative form of equation (\ref{7}):
\begin{equation}\label{11}
\frac{\partial P}{\partial t}-\frac{\partial}{\partial r}
\left(\frac{w^{2}(r)}{2}\frac{\partial P}{\partial r}+
(\theta(r)+(\lambda(t)-w')w)P)\right)+(r+\theta'+\lambda(t)w'-(ww')')P= 0.
\end{equation}

\section{  Interest rate discretization }

Let the interest rate interval $I=(0,R)$ be divided into $N$ sub-intervals
\[
I_{i}:=(r_{i}, r_{i+1}), \; \; i=0,1, \dots, N-1,
\]
with the grid ${\overline{w}}_{h}= w_{h} \cup \{r_{0} \} \cup \{r_{N} \},$
${\overline{w}}_{h} = \{r_{i}, \; i=0, 1, \dots, N, \; 0=r_{0}< r_{1} < \dots < r_{N-1} < r_{N} = R \}$.
For each $i=0,1, \dots , N-1$
we put $h_{i}=r_{i+1}-r_{i}$ and
$h=\max_{0\leq i \leq N-1}h_{i}$.
We also let
$r_{i-1/2}=(r_{i-1}+r_{i})/2$ and $r_{i+1/2}=(r_{i}+r_{i+1})/2$ for each
$i=1, 2, \dots, N-1 $. These mid-points form a {\emph second partition} grid
${\overline{w}}_{{\widetilde{h}}}= w_{{\widetilde{h}}} \cup\{ r_{-1/2}\} \cup \{r_{N+1/2}\}$
of $[0,R ]$ if we define $r_{-1/2}=r_{0}$ and
$ r_{N+1/2}= r_{N}$.
Let $\hbar_{i}=r_{i+1/2}-r_{i-1/2}=0.5(h_{i}+h_{i-1}),$
$\hbar_{0}=h_{1}/{2}=r_{1/2},\;\hbar_{N}=R- r_{N-1/2}.$

According to the assumptions (\ref{9}), \eqref{10}, at the construction of the
finite volume approximation several cases must be considered.

{\bf{Case 1.}} {\emph {We consider equation} (\ref{11}) {\emph{with coefficients} (\ref{9}), \eqref{10a}.}}
Now (\ref{11}) takes the form
\begin{align}\label{12}
\frac{\partial P}{\partial t}-\frac{\partial}{\partial r}
\left[r^{2}(R-r)^{2}\frac{w_{0}^{2}(r)}{2}\frac{\partial P}{\partial r}
+r(R-r)\left(\theta_{0}(r)+(\lambda(t)- w') w_{0}(r)\right)P\right]\nonumber \\
+(r+\theta'+\lambda(t)w'-(ww')')P=0.
\end{align}
Integrating  (\ref{12}) over the interval $(r_{i-1/2}, r_{i+1/2})$ we have
\begin{equation}\label{13}
\int_{r_{i-1/2}}^{r_{i+1/2}}\frac{\partial P}{\partial t}dr-\left[
r(R-r)\left(\frac{w_{0}^{2}(r)}{2}r (R-r)\frac{\partial P}{\partial r}+ bP\right)
\right]_{r_{i-1/2}}^{r_{i+1/2}}+Q_{i}=0,
\end{equation}
\[
Q_{i}=\int_{r_{i-1/2}}^{r_{i+1/2}}(r+\theta'+\lambda(t)w'-(ww')')P dr=0,
\]
for $i=1, 2, \ldots, N-1$, where we denoted
\begin{equation}\label{14}
b=b(r,t)=\theta_{0}(r)+(\lambda(t)-w'(r))w_{0}(r).
\end{equation}
Applying the mid-point qudrature rule to the first and the last terms in (\ref{13}) we obtain
\begin{equation}\label{15}
\frac{\partial P_{i}}{\partial t}\hbar_{i}-\left[r_{i+1/2}(R-r_{i+1/2})\rho(P) \vert_{r_{i+1/2}}-
r_{i-1/2}(R-r_{i-1/2})\rho(P)\vert_{r_{i-1/2}}\right]+Q_{i}^{h}P_{i}=0,
\end{equation}

\[
Q_{i}^{h}=r_{i}\hbar_{i}+\theta_{i+1/2}-\theta_{i-1/2}+\lambda(t)(w_{i+1/2}-w_{i-1/2})
-(ww')_{i+1/2}+(ww')_{i-1/2}
\]
for $i=1, 2, \dots, N-1 $, {\emph {$P_{i}$ denotes the nodal approximation to $P(r_{i}, t)$}}
to be determined and $\rho (P)$ is the flux associated with $P$ and denoted by
\begin{equation}\label{16}
\rho(P):=ar(R-r)\frac{\partial P}{\partial r}+bP, \; \;
a=a(r)=\frac{w_{0}^{2}(r)}{2}.
\end{equation}

The discussion is divided into {\emph{three sub-cases}}.

{\bf{Case 1.1.}} {\emph {Approximation of $\rho$ at $r_{i+1/2}$ for $1 \leq i \leq N-2.$}}

Let us consider the following two-point boundary value problem for $r \in I_{i}$:
\begin{subequations}\label{17}
\begin{align}
(a_{i+1/2}r (R-r)v'+b_{i+1/2}v)'=0, \; \;\label{17a}\\
v(r_{i})=P_{i},\;\; \;\;v(r_{i+1})=P_{i+1},\label{17b}
\end {align}
\end {subequations}
where $a_{i+1/2}=a(r_{i+1/2}),\; b_{i+1/2}= b(r_{i+1/2}, t)$. Integrating \eqref{17a}
yields the first order linear equation
\begin{equation}\label{18}
\rho_{i}(v):=a_{i+1/2}r(R-r)v'+b_{i+1/2}v=C_{1},
\end {equation}
where $C_{1}$ denotes an additive constant (depending on $t$). The analytic solution of this linear equation is
\begin{equation}\label{19}
v(r)=\frac{C_{1}}{b_{i+1/2}}+C_{2}\left(\frac{r}{R-r}\right)^{-\frac{b_{i+1/2}}{R a_{i+1/2}}},
\end{equation}
where $C_{2}$ is an additive constant. Note that in this reasoning we assume that
$b_{i+1/2}\not=0$. But as will be seen below, the restriction can be lifted as it
is limiting case of the above when $b_{i+1/2} \to 0$. Applying the boundary condition \eqref{17b} we obtain
\begin{equation}\label{20}
P_{i}=\frac{C_{1}}{b_{i+1/2}}+C_{2}\left(\frac{r_{i}}{R- r_{i}}\right)^{\frac{-\alpha_{i}}{R}},
\;\;\;\;
P_{i+1} =\frac{C_{1}}{b_{i+1/2}}+C_{2}\left(\frac{r_{i+1}}{R-r_{i+1}}\right)^
{\frac{-\alpha_{i}}{R}},
\end{equation}
where $\alpha_{i}=b_{i+1/2}/a_{i+1/2}$. Solving this linear system gives
\begin{equation}\label{21}
\rho_{i}(P)=C_{1}=b_{i+1/2}\frac{ \left(\frac{r_{ i+1}}{R - r_{i+1}}
\right)^{\frac{\alpha_{i}}{R}} P_{i+1}-\left(\frac{r_{i}}{R-r_{i}}
\right)^{\frac{\alpha_{i}}{R}}P_{i}}{\left(\frac{r_{i+1}}{R-r_{i+1}}
\right)^{\frac{\alpha_{i}}{R}}-\left(\frac{r_{i}}{R-r_{i}}
\right)^{\frac{\alpha_{i}}{R}}}
\end{equation}
for $ i=1, \dots, N-2$.

This gives a representation for the flux on the right-hand side of (\ref{18}).
Note that (\ref{21}) also holds when $\alpha_{i}\to 0$ . This is because
\begin{align}\label{22}
\lim_{\alpha_{i} \to 0}\frac{\left(\frac{r_{ i+1}}{R-r_{i+1}}\right)^
{\frac{\alpha_{i}}{R}}-\left(\frac{r_{i}}{R-r_{i}}\right)^{\frac{\alpha_{i}}{R}}}
{b_{i+1/2}}=\frac{1}{a_{i+1/2}}\lim_{\alpha_{i} \to 0}\frac{\left(\frac{r_{i+1}}{R-r_{i+1}}
\right)^{\frac{\alpha_{i}}{R}}-\left(\frac{r_{i}}{R-r_{i}}\right)^
{\frac{\alpha_{i}}{R}}}{\alpha_{i}}\nonumber\\
= \frac{1}{Ra_{i+1/2}}\left(\ln\frac{r_{i+1}}{R-r_{i+1}}-\ln\frac{r_{ i }}
{R-r_{i}}\right)=\frac{1}{Ra_{i+1/2}}\ln\left(\frac{r_{i+1}}{r_{i}}
\frac{R-r_{i}}{R-r_{i+1}}\right)>0
\end{align}
since $r_{i}<r_{i+1}$ and $a_{i+1/2}>0.$ Thus, $\rho_{i}(P)$ in (\ref{21})
provides an approximation to the flux $\rho_{i}(P)$ at $r_{i+1/2}$.

{\bf{Case 1.2.}} {\emph{Approximation of $\rho$ at $r_{1/2}$.}}

Now, we write the flux in the form
\[
\rho(P):=ar\frac{\partial P}{\partial r}+ bP, \;\;
a=a(r)=\frac{w_{0}^{2}(r)}{2}(R-r).
\]

Note that the analysis in Case 1.1 does not apply to approximation of the
flux because  \eqref{17a} is degenerate. This can be seen from
expression (\ref{19}). When $\alpha_{0}>0$, we have to chose
$ C_{2} = 0 $ as, otherwise, $v$ blows up as  $r \to 0$ .
However, the resulting solution $v=C_{1}/b_{1/2}$ can never satisfy both of conditions in
\eqref{17b}. To solve this difficulty, following [15], we
will reconsider \eqref{17a}, \eqref{17b} with an extra degree of freedom in the following form:
\begin{align*}
(a_{1/2}rv'+ b_{1/2}v)'=C_{2},\;\mbox{in} \;(0,r_{1}),\\
v(0)=P_{0}, \;\;\;\; v(r_{1})=P_{1},
\end{align*}
where $C_{2}$ is an unknown constant to be determined. Integrating the differential
equation once we have
\[
a_{1/2}rv'+b_{1/2} v=C_{2}r+C_{3}.
\]
Using the condition $v(0)=P_{0}$ we have $C_{3}=b_{1/2}P_{0},$ and so the above equation becomes
\begin{equation}\label{23}
\rho_{0}(v):=a_{1/2}rv'+b_{1/2}v=C_{2}r+b_{1/2}P_{0}.
\end{equation}
Solving this problem analytically gives
\begin{equation}\label{24}
v(r)=\left\{\begin{array}{ll}
P_{0}+\frac{C_{2}r}{a_{1/2}+ b_{1/2}}+C_{4}r^{-\alpha_{0}},&\alpha_{0}\not=- 1,\\
\\
P_{0}+\frac{C_{2}}{a_{1/2}}r\ln r+C_{4}r,&\alpha_{0}=- 1,
\end{array}
\right.
\end{equation}
where $\alpha_{0}= b_{1/2}/a_{1/2}$ as defined Case 1.1 and $C_{4}$ is an additive
constant (depending on $t$).

To determine the constant $C_{2}$ and $C_{4}$, we first consider the case $\alpha_{0} \not=- 1$.
When $\alpha_{0}\geq 0,\;v(0)=P_{0}$ implies that $C_{4}=0$. If $\alpha_{0}< 0$, $C_{4}$ is arbitrary,
so we also choose $C_{4}=0$. Using $v(r_{1})=P_{1}$ we obtain
$C_{2}=\frac{1}{r_{1}}\left(a_{1/2}+b_{1/2}\right)(P_{1}-P_{0})$.

When $\alpha_{0}=- 1$, from (\ref{24}) we see that $v(0)=P_{0}$  is satisfied for
any $C_{2}$ and $C_{4}$. Therefore, solutions with such $C_{2}$ and $C_{4}$ are
not unique. We choose $C_{2}=0$, and $v(r_{1})=P_{1}$ and then
$C_{4}=(P_{1}-P_{0})/r_{1}$. Therefore, from (\ref{23}) we have that
\begin{equation}\label{25}
\rho_{0}(P):=(a_{1/2}rv'+b_{1/2}v)_{r_{1/2}}=\frac{1}{2}[(a_{1/2}+
b_{1/2})P_{1}-(a_{1/2}-b_{1/2})P_{0}]
\end{equation}
for both $\alpha_{0}=-1$ and $\alpha_{0}\not= -1$. Furthermore, (\ref{24}) reduces to
\begin{equation}\label{26}
v=P_{0}+(P_{1}-P_{0})r/r_{1},\; r\in [0, r_{1}].
\end{equation}

{\bf{Case 1.3.}} {\emph{Approximation of $\rho$ at $r_{N-1/2}$}}.

We write the flux in the form
\[
\rho(P):=a(R-r)\frac{\partial P}{\partial r}+bP,\;\;\;
a=a(r)=\frac{w_{0}^{2}(r)}{2}r.
\]

The situation is symmetric to this of Case 1.2. We consider the auxiliary problem:
\begin{align*}
(a_{N-1/2}(R-r)v'+b_{N-1/2}v)'=C_{2},\;\mbox{in}\;(r_{N-1}, R),\\
v(r_{N-1})=P_{N-1},\;\;\;v(R)=P_{N},
\end{align*}
where $C_{2}$ is an unknown constant to be determined. Integrating the
differential equation once we have
\[
a_{N-1/2}(R-r)v'+b_{N-1/2}v=C_{2}r+C_{3}.
\]
Using the condition $v(R)=P_{N}$ we have $b_{N-1/2}P_{N}=C_{2}R+C_{3}$, and so the last equation becomes
\begin{equation}\label{27}
\rho_{N-1}(v):=a_{N-1/2}(R-r)v'+b_{N-1/2}v=- C_{2}(R-r)+b_{N-1/2}P_{N}.
\end{equation}

Solving this problem analytically gives
\begin{equation}\label{28}
v(r) =\left\{
\begin{array}{ll}
P_{N} +C_{2}\frac{(R-r)}{(1-\alpha_{N-1})a_{N-1/2}}
+C_{4}(R- r)^{\alpha_{N-1}},& \alpha_{N-1}\not = 1, \\
\\
P_{N}+\frac{C_{2}(R-r)}{a_{N-1/2}}\ln(R-r)+C_{4}(R-r),& \alpha_{N-1}=1,
\end{array}\right.
\end{equation}
where $\alpha_{ N- 1}=b_{N-1/2}/a_{N-1/2}$ as defined before and $C_{4}$
is an additive constant (dependent on $t$).

To determine the constants $C_{2}$ and $C_{4}$, we first consider the case
when $\alpha_{N-1}\not=1$. When $\alpha_{N-1} < 0,\;v(R)=P_{N}$ implies $C_{4}=0$.
If $\alpha_{N-1} \geq 0,\;C_{4}$ is arbitrary, so we also choose $C_{4}=0$.
Using $v(r_{N-1})=P_{N-1}$ in (\ref{28}) we obtain
$C_{2} =(P_{N-1}-P_{N})(a_{N-1/2}-b_{N-1/2})/(R-r_{N-1})$.

When $\alpha_{N-1}= 1$, from (\ref{28}) we see that
$v(R)=P_{N}$ is satisfied for any $C_{2}$ and $C_{4}$. We choose $C_{2}=0$, and
$v(r_{N-1})=P_{N-1}$ in (\ref{28}) gives
$C_{4}=(P_{N-1}-P_{N})/(R-r_{N-1})$. Therefore, from (\ref{27}) we have
\begin{equation}\label{29}
\rho_{N-1}=\frac{1}{2}\left[(a_{N-1/2}+b_{N-1/2})P_{N}-
(a_{N-1/2}-b_{N-1/2})P_{N-1}\right].
\end{equation}

{\bf{Case 2.}} \emph{Now we consider equation} (\ref{11}) \emph{with coefficients} (\ref{9}), \eqref{10b}.

Following the line in Case 1, we have
\begin{equation}\label{30}
\int_{r_{i-1/2}}^{r_{i+1/2}}\frac{\partial P}{\partial t}dr-
\left[r(r\frac{w_{0}^{2}(r)}{2}(R-r)^{2}\frac{\partial P}{\partial r}+bP)
\right]_{r-1/2}^{r_{i+1/2}}+Q_{i}=0
\end{equation}
for $i=1, 2, \dots, N-1$, where
\[
b(r,t)=\theta_{0}(r)+(\lambda(t)-w'(r))(R-r)w_{0}(r).
\]

{\bf{Case 2.1.}} {\emph{Approximation of $\rho$ at $r_{i+1/2}$ for $1\leq i\leq N-2$.}}

Applying the mid-point quadrature rule to the first and
third terms in (\ref{30}) we find
\[
\frac{\partial P_{i}}{\partial t}\hbar_{i}-
[r_{i+1/2}\rho(P)\vert_{r_{i+1/2}}-r_{i-1/2}\rho(P)\vert_{r_{i-1/2}}]
+Q_{i}^{h}=0
\]
for $i=1, 2, \dots, N-1$, where
\[
\rho(P):=ar( R-r)\frac{\partial P}{\partial r}+bP,\;\;
a=a(r)=\frac{w_{0}^{2}(r)}{2}(R-r).
\]
Further, one can obtain a formula in the form (\ref{21}).

{\bf {Case 2.2.}} {\emph{Approximation of $\rho$ at $r_{1/2}.$}} Now we proceed as in Case 1.2, but
\[
\rho(P)=ar\frac{\partial P}{\partial r}+ bP, \; \;
a=a(r)=\frac{w_{0}^{2}(r)}{2}(R-r)^{2}.
\]

{\bf {Case 2.3}} {\emph{ Approximation of $\rho$ at $r_{N-1/2}$.}}
In this case
\[
\rho(P)=a(R-r)\frac{\partial P}{\partial r}+ bP, \; \;
a=a(r)=\frac{w_{0}^{2}(r)}{2}r(R-r).
\]

{\bf{Case 3.}} \emph{Here we consider equation} (\ref{11}) \emph {with coefficients} (\ref{9}), \eqref{10c}.
In this case the construction is symmetric to this in Case 2 and
we will only present the results.
\[
\int_{r_{i-1/2}}^{r_{i+1/2}}
\frac{\partial P}{\partial t}dr-\left[(R-r)(r^{2}\frac{w_{0}^{2}(r)}
{2}(R-r)\frac{\partial P}{\partial r}+bP)\right]_{r-1/2}^{r_{i+1/2}}+Q_{i}= 0
\]
for $i=1, 2, \dots, N-1$, where
\[
b(r,t)=\theta_{0}(r)+(\lambda(t)-w'(r))rw_{0}(r).
\]

{\bf{Case 3.1.}} {\emph{Approximation of $\rho$ at $r_{i+1/2}$ for $1<i\leq N-2$.}} Now we take
\[
\rho(P)=ar(R-r)\frac{\partial P}{\partial r}+ bP,\; \;
a=a(r)=r\frac{w_{0}^{2}(r)}{2}.
\]

{\bf{Case 3.2.}} {\emph{Approximation of $\rho$ at $r_{1/2}$.}}
In this subcase
\[
\rho(P)=ar\frac{\partial P}{\partial r}+bP,\;\;
a=a(r)=r(R-r)\frac{w_{0}^{2}(r)}{2}.
\]

{\bf{Case 3.3.}} {\emph{Approximation of $\rho$ at $r_{ N-1/2}$.}}
Now we proceed as in Case 1.3 but
\[
\rho(P)=a(R-r)\frac{\partial P}{\partial r}+ bP, \; \;
a=a(r)=r^{2}\frac{w_{0}^{2}(r)}{2}.
\]

{\bf{Case 4.}} \emph{ Here we consider equation} (\ref{11}) \emph {with coefficients}  (\ref{9}), \eqref{10d}.
We have
\[
\int_{r_{i-1/2}}^{r_{i+1/2}}dr-\left[r^{2}(R-r)^{2}
\frac{w_{0}^{2}(r)}{2}\frac{\partial P}{\partial r}+bP\right]
_{r-1/2}^{r_{i+1/2}}+Q_{i}= 0
\]
for $ i=1, 2, \dots, N-1$, where
\[
b(r,t)=\theta_{0}(r)+(\lambda(t)-w'(r))w(r).
\]

{\bf{Case 4.1.}} {\emph{Approximation of $\rho$ at $r_{ i+1/2} \; \mbox{for} \;  1 \leq i \leq N-2$.}}
Now we choose
\[
\rho(P)=ar(R-r)\frac{\partial P}{\partial r}+ bP, \; \;
a=a(r)=r(R-r)\frac{w_{0}^{2}(r)}{2}.
\]

{\bf{Case 4.2.}} {\emph{Approximation of $\rho$ at $r_{1/2}$.}}
We take
\[
\rho(P)=ar\frac{\partial P}{\partial r}+ bP, \; \;
a=a(r)=r(R-r)^{2}\frac{w_{0}^{2}(r)}{2}.
\]

{\bf{Case 4.3.}} {\emph{Approximation of $\rho$ at $r_{ N-1/2}$.}}
We choose
\[
\rho(P)=a(R-r)\frac{\partial P}{\partial r}+ bP, \; \;
a=a(r)=r^{2}(R-r)\frac{w_{0}^{2}(r)}{2}.
\]

Finally, using (\ref{21}), (\ref{25}), (\ref{27}) and (\ref{29}),
depending on the value of $i=0, 1, \dots, N-1$ respectively, we define a
global piecewise constant approximation to $\rho(P)$ by $\rho_{h}(P)$ satisfying
\begin{equation}\label{31}
\rho_{h}(P)=\rho_{i}(P)\;\mbox{if}\; x \in I_{i}
\end{equation}
for $i=0, 1, \dots, N- 1$.

Substituting (\ref{21}) or (\ref{25}) or (\ref{27}) or (\ref{29}),
depending on the value of $i$ respectively, into (\ref{15}) we obtain
\begin{equation}\label{32}
\frac{\partial P_{i}}{\partial t}\hbar_{i} -
e_{i,i-1}P_{i-1}+e_{i,i}P_{i}-e_{i,i+1}P_{i+1}=0,
\end{equation}
where
\begin{align*}
e_{1,0}=0.5r_{1/2}(R-r_{1/2})(a_{1/2}-b_{1/2}),\;\;\;
e_{1,1}=0.5r_{1/2}(R-r_{1/2})(a_{1/2}+ b_{1/2})\\
+Q_{1}^{h}+
r_{3/2}(R- r_{3/2})b_{3/2}\frac{ \left(\frac{r_{1}}{R-r_{1}}
\right)^{\frac{\alpha_{1}}{R}}}{\left(\frac{r_{2}}{R-r_{2}}
\right)^{\frac{\alpha_{1}}{R}}-\left(\frac{r_{1}}{R - r_{1}}
\right)^{\frac{\alpha_{1}}{ R}}},\\
e_{1,2}=r_{3/2}(R-r_{3/2})b_{3/2}\frac{\left(\frac{r_{2}}{R-r_{2}}
\right)^{\frac{\alpha_{1}}{R}}}{\left(\frac{r_{2}}{R-r_{2}}
\right)^{\frac{\alpha_{1}}{R}}-
\left(\frac{r_{1}}{R-r_{1}}\right)^{\frac{\alpha_{1}}{R}}},\\
e_{i,i-1}=r_{i-1/2}(R-r_{i-1/2})b_{i-1/2}\frac{ \left(\frac{r_{i-1}}{R-r_{i-1}}
\right)^{\frac{\alpha_{i-1}}{R}}}{\left(\frac{r_{i}}{R-r_{i}}
\right)^{\frac{\alpha_{i-1}}{R}}-\left(\frac{r_{i-1}}{R-r_{i-1}}
\right)^{\frac{\alpha_{i-1}}{R}}},\\
e_{i,i}=Q_{i}^{h}+r_{i+1/2}(R-r_{i+1/2})b_{i+1/2}
\frac{\left(\frac{r_{i}}{R-r_{i}}\right)^{\frac{\alpha_{i}}{R}}}
{\left(\frac{r_{i+1}}{R-r_{i+1}}\right)^{\frac{\alpha_{i}}{R}}-
\left(\frac{r_{i}}{R - r_{i}}\right)^{\frac{\alpha_{i}}{R}}} \\
+r_{i-1/2}(R-r_{i-1/2})b_{i-1/2}\frac{\left(\frac{r_{i}}{R-r_{i}}
\right)^{\frac{\alpha_{i-1}}{R}}}{\left(\frac{r_{i}}{R-r_{i}}
\right)^{\frac{\alpha_{i-1}}{R}}-\left(\frac{r_{i-1}}{R-r_{i-1}}
\right)^{\frac{\alpha_{i-1}}{R}}},\\
e_{i,i+1}=r_{i+1/2}(R-r_{i+1/2})b_{i+1/2}\frac{ \left(\frac{r_{i+1}}{R-r_{i+1}}
\right)^{\frac{\alpha_{i}}{R}}}{\left(\frac{r_{i+1}}{R-r_{i+1}}
\right)^{\frac{\alpha_{i}}{R}}-\left(\frac{r_{i}}{R-r_{i}}\right)^{\frac{\alpha_{i}}{R}}},
\end{align*}
for $i=2, 3, \dots, N-2$;
\begin{align*}
e_{N-1,N-2}=r_{N-3/2}(R-r_{N-3/2})
b_{N-3/2}\frac{\left(\frac{r_{N-2}}{R-r_{N-2}}
\right)^{\frac{\alpha_{N-2}}{R}}}
{\left(\frac{r_{N-1}}{R-r_{N-1}}
\right)^{\frac{\alpha_{N-2}}{R}}-
\left(\frac{r_{N-2}}{R-r_{N-2}}
\right)^{\frac{\alpha_{N-2}}{R}}},\\
e_{N-1,N-1}=r_{N-3/2}(R-r_{N-3/2})b_{N-3/2}
\frac{ \left(\frac{r_{N-1}}{R-r_{N-1}}
\right)^{\frac{\alpha_{N-2}}{R}}}
{\left(\frac{r_{ N-1}}{R-r_{N-1}}
\right)^{\frac{\alpha_{N-2}}{R}}-
\left(\frac{r_{N-2}}{R-r_{N-2}}
\right)^{\frac{\alpha_{N-2}}{R}}}\\
+0.5r_{ N-1/2}(R-r_{N-1/2})(a_{N-1/2}-b_{N-1/2})+ Q_{N-1}^{h},\\
e_{N-1,N}=0.5r_{N-1/2}(R-r_{N-1/2})(a_{N-1/2}+b_{N-1/2}).
\end{align*}

Now we will derive the semi-discrete equations at
$r=0$ and $r=R$. We integrate the equation (\ref{12})
over the interval $(r_{-1/2}, r_{1/2})=(r_{0}, r_{1/2})=(0, r_{1/2})$ to get
\[
\int_{0}^{r_{1/2}}\frac{\partial P}{\partial t}dr
-r_{1/2}(R-r_{1/2})\rho(P)\vert_{r_{1/2}}
+\int_{0}^{r_{1/2}}(r+\theta'+\lambda(t)w'-(ww')')Pdr=0.
\]
Using (\ref{25}) we obtain
\[
\frac{\partial P_{0}}{\partial t} \frac{h_{0}}{2}
-\frac{h_{0}}{4}\left(R-\frac{h_{0}}{2}\right)
[(a_{1/2}+b_{1/2})P_{1}-(a_{1/2}-b_{1/2})P_{0}]
+Q_{0}^{h}P_{0}=0,
\]
where
\[
Q_{h}^{0}=\frac{h_{0}^{2}}{8}
+\theta(r_{1/2})+\lambda(t)w(r_{1/2})-ww' \vert_{r_{1/2}}.
\]
Therefore, at $r=0$ we have:
\begin{equation}\label{33}
\frac{\partial P_{0}}{\partial t}\frac{h_{0}}{2}+
e_{0,0}P_{0}-e_{0,1}P_{1}=0,
\end{equation}
\[
e_{0,0}=\frac{h_{0}}{4}\left(R-\frac{h_{0}}{2}\right)
(a_{1/2}-b_{1/2})+Q_{0}^{h}, \; \; \;
e_{0,1}=\frac{h_{0}}{4}\left(R-\frac{h_{0}}{2}\right)(a_{1/2}+b_{1/2}).
\]

Next, in a similar way (now integrating (\ref{12}) over
$(r_{N-1/2}, r_{N+1/2})=(r_{N-1/2},\\ r_{N})=
\left(R-\frac{h_{N-1}}{2}, R\right)$
and using (\ref{29})), we derive the semi-discrete equation at
$x_{N}=R$:
\begin{align*}
\frac{\partial P_{N}}{\partial t}\frac{h_{N-1}}{2}+
\frac{h_{N-1}}{4}\left(R-\frac{h_{N-1}}{2}\right)[(a_{N-1/2}+b_{N-1/2})P_{N}-\\
(a_{N-1/2}-b_{N-1/2})P_{N-1}]+Q_{N}^{h}P_{N}=0,
\end{align*}
where
\[
Q_{N}^{h} =\frac{h_{N-1}}{4}
\left(2R-\frac{h_{N-1}}{2}\right)-\theta
\vert_{R-\frac{h_{N-1}}{2}}-
\lambda(t)w\vert_{R-\frac{h_{N-1}}{2}}
+(ww')\vert_{R-\frac{h_{N-1}}{2}}.
\]
Therefore, at $r=R$ we have
\begin{equation}\label{34}
\frac{\partial P_{N}}{\partial t}\frac{h_{N-1}}{2}-
e_{N,N-1}P_{N-1}+e_{N,N}P_{N}=0,
\end{equation}
where
\[
e_{N,N-1}=\frac{h_{N-1}}{4}\left(R-\frac{h_{N-1}}{2}\right)(a_{N-1/2}-b_{N-1/2}),
\]
\[
e_{N,N} =\frac{h_{N-1}}{4}\left(R-\frac{h_{N-1}}{2}\right)(a_{N-1/2}+b_{N-1/2})+Q_{N}^{h}.
\]

We now discuss the accuracy of the interest rate discretization of the system
(\ref{32}), (\ref{33}), (\ref{34}). Let
${\bf{E}}_{i},\;i=0, 1, \dots,N$ be $N+1$ row vectors with dimension $N+1$ defined by
\begin{align*}
{\bf{E}}_{0}(t)=(e_{0,0}(t), -e_{0,1}(t), 0, \dots, 0),\;\;
{\bf{E}}_{N}(t)=(0, \dots, -e_{N,N-1}(t), e_{N,N}(t)),\\
{\bf{E}}_{i}(t)=(0, \dots, - e_{i,i-1}(t), e_{i,i} (t), - e_{i,i+1}(t),0, \dots, 0),
\; \; i=1, 2, \dots, N-1.
\end{align*}
Obviously, introducing the vector ${\bf{P}}=(P_{0}(t), P_{1}(t), \dots, P_{N}(t))^{T}$
and using ${\bf{E}}_{i}$, the equations (\ref{32}), (\ref{33}), (\ref{34}) can be written as
\begin{equation}\label{35}
\frac{dP_{i}(t)}{dt}\hbar_{i}-{\bf{E}}_{i}(t){\bf{P}}(t)=0,
\end{equation}
for $i=0, 1, \dots, N$. This is a first-order linear ODEs system.

To estimate the accuracy of the interest rate discretization, we
will follow \cite{W}. First, we define a space $S_{h}$ of functions $\phi_{i}$ associated with
$r_{i}$ in the following way. On the interval $(r_{i}, r_{i+1})$ we choose $\phi_{i}$
so that it satisfies \eqref{17a} with $\phi_{i}(r_{i})=1$ and $\phi_{i}(r_{i+1})=0$.
Naturally, the solution to this two-point boundary value problem is given in (\ref{19})
where $C_{1}$ and $C_{2}$ are determined by (\ref{20}) with $P_{i}=1$ and $P_{i+1}=0$. Similarly we
define $\phi_{i}(r)$ on the interval $(r_{i-1}, r_{i})$ so that $\phi_{i}(r_{i-1})=0$ and
$\phi_{i}(r_{i})=1$. Combining these two solutions and extending the function $\phi_{i}(r)$ as zero
to the rest of the interval $(0,R)$ we have for $i=1,\dots,N-1$

$$
\phi_{i} (r) =
                   \left\{
                     \begin{array}{ll}
\displaystyle                         \frac{ \left(
                                  \frac{   R } {   r_{i-1}   }
                             -1    \right)^{   \frac{ \alpha_{i-1}}{ R  }  } -
                                \left(
                                  \frac{   R } {  r  } -1
                                \right)^{   \frac{ \alpha_{i-1}}{ R  }  } }
                                {\left(
                                  \frac{  R } {   r_{i-1}   }
                              -1  \right)^{   \frac{ \alpha_{i-1}}{ R  }  }  -
                                \left(
                                  \frac{   R } {   r_{i}   }
                              -1  \right)^{   \frac{ \alpha_{i-1}}{ R  }  } }  ,
                       r \in ( r_{i-1},r_{i} ) ,
                        \\
\displaystyle                         \frac{ \left(
                                  \frac{   R } {   r_{i+1}   }
                             -1    \right)^{   \frac{ \alpha_{i}  }{ R  }  } -
                                \left(\frac{   R } {  r  } -1
                                \right)^{   \frac{ \alpha_{i}   }{ R  }  } }
                                {\left(
                                  \frac{  R } {   r_{i+1}   }
                              -1  \right)^{   \frac{ \alpha_{i}  }{ R  }  }  -
                                 \left(\frac{   R } {   r_{i}   }
                              -1  \right)^{   \frac{ \alpha_{i}  }{ R  }  } },
                       r \in ( r_{i},r_{i+1} ), \\
                       0, \mbox {otherwise}.
                      \end{array}
                      \right.
                      $$

In a similar way, on the intervals $(0, r_{1})$ and $(r_{N-1}, R)$
we define the linear functions
\[
\phi_{0}(r)=\left\{
\begin{array}{ll}
1-\frac{r}{r_{1}},\;\; r \in (0, r_{1}),\\
\\
0,\mbox {otherwise},
\end{array}
\right.
\; \;\; \;\; \;
\phi_{N}(r)=\left\{
\begin{array}{ll}
\frac{r-r_{N-1}}{R-r_{N-1}},\;\; r \in (r_{N-1}, R), \\
\\
0, \mbox {otherwise}.
\end{array}
\right.
\]
The following assertion is an analogue of Lemma 4.2 in \cite{W}.

\begin{lem}
Let $v$ be a sufficiently smooth function and $v_{I}$ be the $S_{h}$-interpolant of $v$.
Then
\[
\Vert \rho(v)-\rho_{h}(v_{I})\Vert_{\infty,I_{i} }\leq
C(\Vert \rho'(v) \Vert_{\infty,I_{i}}+\Vert b'\Vert_{\infty, I_{i}}
\Vert v \Vert_{\infty, I_{i}})h_{i},
\]
$i=0, 1, \dots, N-1$ where $\rho$ and $\rho_{h}$ are the fluxes defined
in (\ref{16}) and (\ref{31}), respectively and $C$ is a positive constant independent of
$h_{i}$ and $v$.
\end{lem}

Summarizing the constructions in all Cases 1-4 and using Lemma 1,
the following result has been established.

\begin{thm}
The semidiscretization (\ref{35}) is consistent with equation (\ref{7}) and the
truncation error is of order $O(h)$.
\end{thm}

\section{Full discretization}

To discretize the system (\ref{35}) we introduce the time mesh:
\[
{\overline {w}}_{\tau}= w_{\tau}\cup \{0\} \cup \{T\}, \; \;
{\overline {w}}_{\tau}=\{0=t_{0}< t_{1}< \dots<t_{M}= T\}.
\]
For each $j=0, 1, \dots, M-1$ we put $\tau_{j}=t_{j+1}-t_{j}$ and
$\tau=\max_{0 \leq j \leq M-1} \tau_{j}$. Then, we apply the two-level time-stepping method with
splitting parameter $\xi \in [0, 1]$ to (\ref{35}) and yield
\[
\frac{{P}_{i}^{j+1}-{P}_{i}^{j}}{\tau_{j}}\hbar_{i}+\xi{\bf{E}}_{i}^{j+1}{\bf{P}}^{ j+1}
+(1-\xi){\bf{E}}_{i}^{j}{\bf{P}}^{j} = 0
\]
for $j=0, 1, \dots, M-1$. This linear system can be rewritten as
\begin{equation}\label{36}
(\xi{\bf{E}}^{j+1}+{\bf{G}}^{j}){\bf{P}}^{j+1}=
[{\bf{G}}^{j}-(1-\xi){\bf{E}}^{j}]{\bf{P}}^{j}
\end{equation}
for $j=0, 1, \dots, M-1$, where
\[
{\bf{G}}^{j}=diag\left(\frac {h_{0}}{2\tau_{j}},\frac{\hbar_{1}}{\tau_{j}}, \dots,
\frac{\hbar_{N-1}}{\tau_{j}},\frac{h_{N-1}}{2\tau_{j}}
\right)
\]
is $(N+1)\times (N+1)$ diagonal matrix. When $\xi= 1/2$, the time stepping scheme
becomes Crank-Nicholson scheme and when $\theta=1$ it is the backward Euler scheme.
Both of these schemes are unconditionally stable, and they are of second and first order accuracy \cite{W}.

We now show that, when $\tau_{j}$ is sufficiently small, the system matrix of (\ref{36}) is an $M$-matrix.

\begin{thm}
For any given $j=1, 2, \dots, M-1$, if $\tau_{j}$ is sufficiently small, the system matrix of (\ref{36}) is an $M$-matrix.
\end{thm}

\begin{pf}
We will proceed as follows. Using the definition of ${\bf{E}}_{i}^{j+1}$, $i=0,1, \dots , N-1, N $ we will write down the scalar form of (\ref{36}):

$$
B_{0} P_{0}^{j+1} + C_{0}P_{1}^{j+1} = F_{0}
$$

$$
A_{1} P_{0}^{j+1} +
B_{1} P_{1}^{j+1} +
C_{1} P_{2}^{j+1} = F_{1}
$$

$$
A_{2} P_{1}^{j+1} +
B_{2} P_{2}^{j+1} +
C_{2} P_{3}^{j+1} = F_{2}
$$

$$
...............................
$$

$$
A_{i} P_{i-1}^{j+1} +
B_{i} P_{i}^{j+1} +
C_{i}  P_{i+1}^{j+1} = F_{i}
$$

$$
...............................
$$

$$
A_{N} P_{N-1}^{j+1} +
B_{N} P_{N}^{j+1}
= F_{N},
$$

where

$$
B_{0} = \frac{h_{0}}{2 \tau_{j}} + \xi e_{0,0}
, \; \;
C_{0} =
-\xi e_{0,1}
, \; \;
$$

$$
A_{1} =
- \xi e_{1,0}
, \; \;
B_{1} =
 \frac{  \hbar_{1} }{   \tau_{j}  }
 + \xi  e_{1,1}
, \; \;
C_{1} =
-\xi e_{1,2},
$$

$$
A_{i} =
- \xi e_{i,i-1}
, \; \;
B_{i} =
 \frac{  \hbar_{i} }{   \tau_{j}  }
 + \xi  e_{i,i}
, \; \;
C_{i} =
- \xi e_{i,i+1}
, \; \;
i=2,3,\ldots,N-1,
$$


$$
A_{N} =
  - \xi e_{N,N-1}
, \; \;
B_{N} =
 \frac{  h_{N-1} }{  2 \tau_{j}  }
 + \xi e_{N,N}
, \; \;
$$


$$
F_{0} =
\left(
 \frac{  h_{0} }{  2 \tau_{j}  }
 -  ( 1 - \xi ) e_{0,0}
 \right)
  P_{0}^{j} +
  ( 1 - \xi )
  e_{0,1}  P_{1}^{j},
  $$

  $$
  F_{1}
 =
  ( 1 - \xi ) e_{1,0} P_{0}^{j} +
      \left(
 \frac{  \hbar_{1} }{ \tau_{j} }
 -  ( 1 - \xi ) e_{1,1}
 \right)
  P_{1}^{j} +
  ( 1 - \xi )
  e_{1,2}  P_{2}^{j},
  $$

  $$
  F_{i} =
( 1 - \xi ) e_{i,i-1} P_{i-1}^{j} +
     \left(
 \frac{  \hbar_{i} }{   \tau_{j}  }
 -  ( 1 - \xi ) e_{i,i}
 \right)
  P_{i}^{j} +
  ( 1 - \xi )
  e_{i,i+1}  P_{i+1}^{j},
  $$

  $$
  F_{N} =
( 1 - \xi ) e_{N,N-1} P_{N-1}^{j} +
      \left(
 \frac{  h_{N-1} }{  2 \tau_{j}  }
 -  ( 1 - \xi ) e_{N,N}
 \right)
  P_{N}^{j}.
  $$

Let us first investigate the off-diagonal entries of the system matrix $A_{i} = - \xi e_{i,i-1}$ and $C_{i} = - \xi e_{i,i+1}$.
From the formulas for $ e_{i,l} $ from the above we have  $e_{i,l}> 0$,  $ i, l = 1, 2, \dots, N-1, i \ne l.$ That is because
\[
                         b_{   i+1/2 }
                         \frac { \left(
                                  \frac{   r_{ i+1 } } {  R - r_{i+1}   }
                                \right)^{   \frac{ \alpha_{i}  }{ R  }  }}
                                {
                                \left(
                                  \frac{   r_{ i+1 } } {  R - r_{i+1}   }
                                \right)^{   \frac{ \alpha_{i}  }{ R  }  }  -                                                             \left(
                                  \frac{   r_{ i } } {  R - r_{i}   }
                                \right)^{   \frac{ \alpha_{i}  }{ R  }  }  }=
                                a_{i+1/2}\alpha_{i}\frac { \left(
                                  \frac{   r_{ i+1 } } {  R - r_{i+1}   }
                                \right)^{   \frac{ \alpha_{i}  }{ R  }  }}
                                {
                                \left(
                                  \frac{   r_{ i+1 } } {  R - r_{i+1}   }
                                \right)^{   \frac{ \alpha_{i}  }{ R  }  }  -                                                             \left(
                                  \frac{   r_{ i } } {  R - r_{i}   }
                                \right)^{   \frac{ \alpha_{i}  }{ R  }  }  }
                               =
\]
$$
= a_{i+1/2}
\frac{  \alpha_{i} }
     { 1 - { \overline { r}}_{  i } ^{   \frac{ \alpha_{i}  }{  R  }    }       } > 0, \; \;
     0 < { \overline {r}}_{i}  =
     \frac{ r_{i}  }{  r_{ i+1 }  } \cdot
     \frac{  R - r_{  i+1 }  }{ R - r_{i}  } < 1
     $$

\noindent for each $ i=1,2, \dots , N-1 $ and each $ b_{ i+1/2  } \not = 0 $.
We have used that
$ 1 - { \overline { r}}_{i}^{    \frac{ \alpha_{i}   }{R}  } $ has just the  sign of $ \alpha_{i} $.
From (22) we have that it is true also for $ b_{ i+1/2  } \to  0$.
Now it is clear that $A_{i} = - \xi e_{i,i-1}$ and $C_{i} = - \xi e_{i,i+1}$ are negative.

We should also note that $B_{i}$ is always positive since $\tau_{j}$ is small.

 The situation is different for $B_{0}$, $C_{0}$, $A_{1}$, $B_{1}$, $C_{1}$
and
$A_{N-1}$, $B_{N-1}$, $C_{N-1}$, $A_{N}$, $B_{N}$.
       From the first three equations we find

       $$
       P_{0}^{j+1} = \frac{F_{0}}{  B_{0}  } - \frac{ C_{0}} {   B_{0} } P_{1}^{j+1} ,
       \; \;
       P_{1}^{j+1} =
       \frac{ \triangle_{1} }{ \triangle }
       -
       \frac{ C_{1} }{ \triangle } P_{2}^{j+1} ,
       $$

       $$
       \triangle = B_{1} -
       \frac{ A_{1} }{ B_{0} } C_{0} , \; \;
       \triangle_{1} = F_{1}
       - \frac{  A_{1} }{  B_{0} } F_{0},
       $$

       $$
      { \widetilde {B}}_{2}  P_{2}^{j+1} + C_{2} P_{3}^{j+1} =
      { \widetilde {F}}_{2},
      $$

      $$
      { \widetilde {B}}_{2}=
      B_{2} -
      \frac{  A_{2} C_{1} } { \triangle },
     \; \;
      { \widetilde {F}}_{2} =
      F_{2} -
      \frac{ \triangle_{1} }{ \triangle } A_{2}.
      $$

      It is easily to see that when $ \triangle > 0 $ and
      $ \triangle = O
      \left(
       \frac{1}  {   \tau_{j} }
      \right) $ then
      $ B_{2} = O \left( \frac{1} { \tau_{j}   }  \right) $
      for small $ \tau_{j} $.
      Therefore
      $ { \widetilde {B}}_{2} > O \left( \frac{1} { \tau_{j} } \right)$ and $ { \widetilde {B}}_{2} > \vert C_{2} \vert $.

      In a similar way one can eliminate $ P_{N-1}^{j+1}$ and $P_{N}^{j+1} $.
      As a result we obtain a system of linear algebraic equations with unknowns
$ P_{2}^{j+1} , \dots , P_{N-2}^{j+1} $ which matrix is a $M$-matrix.

%

    While $F_{3},...,F_{N-3}$ are non-negative, we have to prove if $\widetilde {F}_{2}$ and $\widetilde {F}_{N-2}$ are also non-negative. From the formulae for $\widetilde {F}_{2}$ it follows that when $\tau_{j}$ is small $\widetilde {F}_{2}$ is non-negative since $F_{2}=O \left( \frac{1} { \tau_{j}   }  \right) $ and ${ \triangle },{ \triangle_{1} }$ are of the same order with respect to $\tau_{j}$. $\widetilde {F}_{N-2}$ is being handled the same way as $\widetilde {F}_{2}$ and also considered non-negative.

    Since the load vector $( \widetilde {F}_{2},F_{3}, \dots ,F_{N-3},\widetilde {F}_{N-2} )$ is non-negative and the corresponding matrix is an M-matrix we can conclude that $ P_{2}^{j+1} , \dots , P_{N-2}^{j+1} $ are non-negative. Finally, using the formulas for $P_{0}^{j+1},P_{1}^{j+1},P_{N-1}^{j+1},P_{N}^{j+1}$ one can easily check that they are non-negative too if $\tau_{j}$ is small.

\end{pf}

\begin{rmk}
    Theorem 3 shows that the fully discretized system (\ref{36}) satisfies
    the discrete maximum principle and because of that fact the above discretization is monotone.
    This guarantees the following:
    for non-negative initial function
    $P_{0} $ the numerical solution $ P_{i}^{j} $, obtained via this method,
    is also non-negative as expected, because the price of the bond is a positive number, see Lemma 1.
\end{rmk}

\section{Numerical Experiments}

Numerical experiments presented in this section illustrate the
properties of the constructed schemes.
In order to investigate numerically the convergence and the
accuracy of the constructed schemes for $\xi=0$, $\xi=1$ and
$\xi=0.5$ we approximately solve the model problem with the known analytical solution
$u(r,t)=\exp(-r-t)$ (exponentially decreasing with respect
to the arguments). We choose this function because its feature is
similar to that of the exact solution to the problem under consideration.
We take $R=1$ and $T=1$. The initial distribution $P_0 (x)$ we compute using this
analytical solution. Let us note, that when we use analytical
solution, in the equation a right hand side arises.

In the tables below are presented the calculated $C$, $L_2$ and $H_1$  mesh norms of the error $z=P-u$ by the formulas
\[
\left\|z \right\|_C =\mathop{\max}\limits_{i,j}\Vert{P_i^j - u_i^j } \Vert
/\max_{i,j} \Vert P_{i}^{j} \Vert,\;\;\;
\left\|z \right\|_{L_2}=\sqrt{\sum\limits_{i = 0}^N
{\sum\limits_{j=0}^M {h\tau \left({P_i^j - u_i^j}\right)}^2}},
\]

\[
\left\|z\right\|_{H_1}=\sqrt{\sum\limits_{i=1}^{N-1}
{\sum\limits_{j = 0}^M {h\tau\left[{\left({P_i^j-u_i^j}
\right)^2+\left({{P'}_i^j- u^j_{\mathop x\limits^\circ,i}}\right)^2} \right]}}}.
\]
Everywhere the calculations are performed with constant time step $\tau = 0.001$.
For the first and the second examples the rate of convergence (RC) is calculated using double mesh principle
\[
    RC=\log_{2}(ER^{N}/ER^{2N}),\;\; ER^{N}=\|P^N-u^N\|,
\]
where $\|.\|$ is the mesh $C$-norm,  $L_{2}$-norm or $H_{1}$-norm, $u^N$ and $P^N$ are respectively the exact solution and the numerical solution computed at the mesh with $N$ subintervals.

\emph{First example.} For the first example coefficients in
equation (\ref{7}) are
\[
\omega(r)=r(R-r),\;\;\; \theta(r)=r(R-r),\;\;\; \lambda(t)=0.25(1+t^2)^{-1}.
\]
That correspond to Case 1. In  Table 1 below are presented the calculated $C$, $L_2$ and $H_1$  mesh norms of the error.
\begin{table}[htbp]
\caption{Crank-Nicholson scheme}
\centering
\begin{tabular} {|c|c|c|c|c|c|c|}
\hline
N&$C$-norm&RC&$L_2$-norm&RC&$H_1$-norm&RC \\
\hline
21&1.481 E-2&-&2.552 E-3&-&2.725 E-2&- \\
\hline
41&7.607 E-3&0.96&9.415 E-4&1.44&1.978 E-2&0.46 \\
\hline
81&3.855 E-3&0.98&3.402 E-4&1.47&1.418 E-2&0.48 \\
\hline
161&1.941 E-3&0.99&1.216 E-4&1.48&1.010 E-2&0.49 \\
\hline
321&9.738 E-4&1.00&4.324 E-5&1.49&7.169 E-3&0.49 \\
\hline
\end{tabular}
\label{tab1}
\end{table}

\emph{Second example.} For the second example coefficients in equation (\ref{7}) are
\[
\omega(r)=r(R-r),\;\;\; \theta(r)=r(R-r)(0.5R-r),\;\;\; \lambda(t)=0.25(1+t^2)^{-1}.
\]
That correspond also to Case 1. In Table 2 below are calculated the mesh $C$, $L_2$ and $H_1$ norms of the error.
\begin{table}[htbp]
\caption{Crank-Nicholson scheme}
\centering
\begin{tabular} {|c|c|c|c|c|c|c|}
\hline
N&$C$-norm&RC&$L_2$-norm&RC&$H_1$-norm&RC \\
\hline
21&1.003 E-2&-&1.482 E-3&-&1.541 E-2&- \\
\hline
41&5.156 E-3&0.96&5.443 E-4&1.44&1.111 E-2&0.46 \\
\hline
81&2.614 E-3&0.98&1.962 E-4&1.47&7.937 E-3&0.48 \\
\hline
161&1.316 E-3&0.99&7.005 E-5&1.48&5.641 E-3&0.49 \\
\hline
321&6.604 E-4&0.99&2.489 E-5&1.49&3.998 E-3&0.49 \\
\hline
\end{tabular}
\label{tab2}
\end{table}

It can be seen from Table 1 and Table 2 that the numerical results are similar.

\emph{Third example.} For this example the coefficients in equation (\ref{7}) are the following :
\[
\omega(r)=r(R-r),\;\;\; \theta(r)=0.5R-r,\;\;\; \lambda(t)=0.25(1+t^2)^{-1},
\]
that correspond to Case 4. Let us note that this case is the most complicated of the four cases discussed in the article with respect to the deriving of the numerical scheme.

In Figure 1 we present the analytical and corresponding approximate solutions.
One can see that the biggest error is \emph{near the ends of the interval, i. e. near to the points of the degeneration}.

\begin{figure}[htbp]
\centering
\includegraphics[width=4in,height=2in]{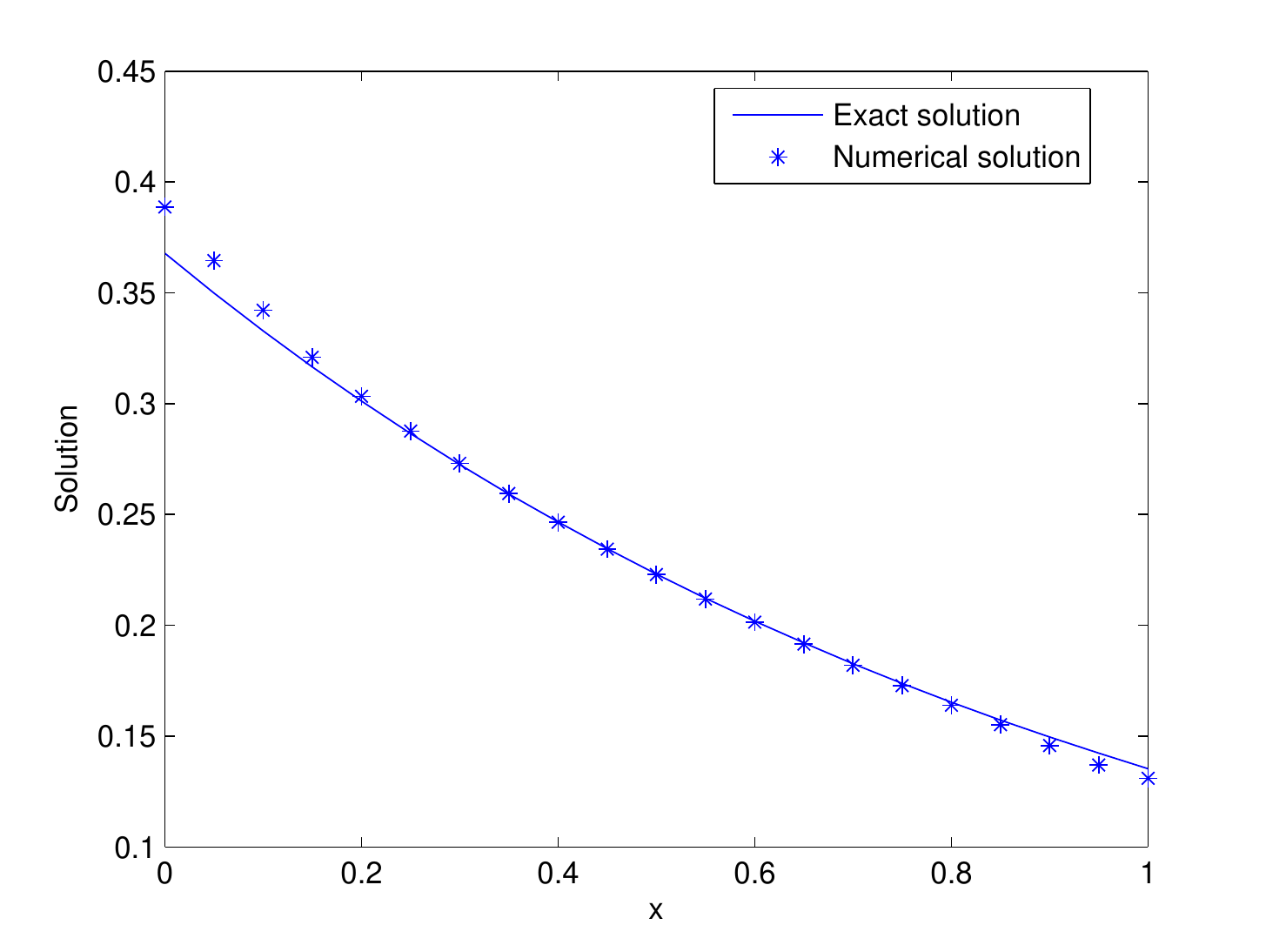}
\caption{Analytical solution $u=exp(-r-t)$, numerical solution for $N=20$, $t=1$.}\label{Fig3}
\end{figure}
In Table 3 are presented the calculated  mesh $C$, $L_2$ and $H_1$ norms of the error for this example.

\begin{table}[htbp]
\caption{Implicit scheme}
\centering
\begin{tabular} {|c|c|c|c|}
\hline
N&$C$-norm&$L_2$-norm&$H_1$-norm \\
\hline
21&2.253 E-2&3.498 E-3&4.078 E-2 \\
\hline
41&8.382 E-3&1.771 E-3&3.561 E-2 \\
\hline
81&4.920 E-3&8.342 E-4&2.728 E-2 \\
\hline
161&2.732 E-3&3.735 E-4&1.965 E-2 \\
\hline
\end{tabular}
\label{tab1}
\end{table}

For this example we used Runge method for practical estimation of the \emph{rate of convergence} $s$
of the considered schemes with respect to the space variable at fixed value of $t$.
In the case when the exact solution $u(x,t)$ of the model problem is known the formula for $s$ is
\[
s=\ln{\left|\frac{u(r)-P_{h}(r)}{u(r)-P_{h/2}(r)}\right|}/\ln{2},
\]
and in the case when the exact solution is not known the formula for $s$ is
\[
s=\ln{\left|\frac{P_{h}(r)-P_{h/2}(r)}{P_{h/2}(r)-P_{h/4}(x)}\right|}/\ln{2}.
\]
In both cases - on two inserted grids (when use the exact solution $u(r,t)$ of model problem) and
on three inserted grids (without exact solution) we get that the rate of convergence is about two,
when the node is not very near to the points of degeneration.

For the problem under consideration we constructed several difference sche\-mes,
well known for non-degenerate parabolic problems [11].
Then, the differential equation (\ref{7}) was approximated, together with the boundary
conditions (\ref{5}), (\ref{6}) and initial condition (\ref{8}).
With respect to the variable $r$ for approximation of the second
derivative is used the usual three-point approximation, and for
the first derivative - central difference. With respect to time a
Crank-Nicolson scheme is constructed. Further this scheme we will
call $B$ Scheme. The scheme we have constructed in this paper for the Case 4 we will call $A$
scheme. From the Table 4 one can see that the scheme A gives more
accurate results near the ends of the interval,
where the degeneration occurs.

\begin{table}[htbp]
\caption{Comparison between Scheme A and Scheme B}
\begin{tabular}{|c|c|c|c|}
\hline
Total Points;Current Point&Time&Scheme A&Scheme B\\
\hline
41&$T$=0.25&&\\
\hline
0&&1.773 E-003&7.874 E-003\\
1&&2.483 E-003&1.216 E-002\\
39&&3.263 E-003&4.157 E-003\\
40&&7.607 E-004&3.071 E-003\\
\hline
81&$T$=0.25&&\\
\hline
0&&3.224 E-004&4.955 E-003\\
1&&8.274 E-006&5.268 E-003\\
79&&1.873 E-003&1.868 E-003\\
80&&8.850 E-006&1.823 E-003\\
\hline
161&$T$=0.25&&\\
\hline
0&&3.405 E-004&4.955 E-003\\
1&&2.897 E-004&5.268 E-003\\
159&&9.900 E-004&1.868 E-003\\
160&&7.775 E-005&1.823 E-003\\
\hline
\end{tabular}
\label{tab3}
\end{table}


\section{Conclusions}

We have studied a degenerate parabolic equation in the zero-coupon bond pricing.
We constructed  and discussed a finite volume difference scheme for
the problem. We have shown that the numerical scheme results a monotone numerical
scheme. The numerical experiments demonstrate the efficiency of our scheme near degeneration.

\vspace{0.5cm}

{ \bf {  Acknowledgements  }}

The first author is supported by the Sofia University Foundation under
Grant No 154/2011. The second author is supported by the Project Bg-Sk-203.

\end{document}